\documentclass[12pt]{article}
\usepackage{amsfonts,amsmath,amsxtra}
\usepackage{amssymb}
\def\hybrid{\topmargin 0pt      \oddsidemargin 0pt
     \headheight 0pt \headsep 0pt
     \textwidth 16.5cm
     \textheight 23cm
     \marginparwidth 0.0in
     \parskip 5pt plus 1pt   \jot = 1.5ex}
\catcode`\@=11
\def\marginnote#1{}
\newcount\hour
\newcount\minute
\newtoks\amorpm
\hour=\time\divide\hour by60
\minute=\time{\multiply\hour by60 \global\advance\minute by-\hour}
\edef\standardtime{{\ifnum\hour<12 \global\amorpm={am}%
     \else\global\amorpm={pm}\advance\hour by-12 \fi
     \ifnum\hour=0 \hour=12 \fi
   \number\hour:\ifnum\minute<10 0\fi\number\minute\the\amorpm}}
\edef\militarytime{\number\hour:\ifnum\minute<10 0\fi\number\minute}
\def\draftlabel#1{{\@bsphack\if@filesw {\let\thepage\relax
\xdef\@gtempa{\write\@auxout{\string
   \newlabel{#1}{{\@currentlabel}{\thepage}}}}}\@gtempa
\if@nobreak \ifvmode\nobreak\fi\fi\fi\@esphack}
     \gdef\@eqnlabel{#1}}
\def\@eqnlabel{}
\def\@vacuum{}
\def\draftmarginnote#1{\marginpar{\raggedright\scriptsize\tt#1}}
\def\draft{\oddsidemargin -0.1truein
     \def\@oddfoot{\sl preliminary draft \hfil
     \rm\thepage\hfil\sl\today\quad\militarytime}
     \let\@evenfoot\@oddfoot \overfullrule 3pt
     \let\label=\draftlabel
     \let\marginnote=\draftmarginnote
\def\@eqnnum{{\rm (\theequation)}
\rlap{\kern\marginparsep\tt\@eqnlabel}%
\global\let\@eqnlabel\@vacuum}  }
\newcommand{\RR}{{\mathbb{R}}}
\newcommand{\CC}{{\mathbb{C}}}

\newcommand{\ZZ}{{\mathbb{Z}}}

\newfont{\Bbbb}{msbm7 scaled 1\@ptsize00}
\newcommand{\zs}{\raise-1pt\hbox{$\mbox{\Bbbb Z}$}}

\font\sevenmsa=msam6 
\newfam\msafam
\textfont\msafam=\sevenmsa
\def\hexnumber@#1{\ifnum#1<10 \number#1\else
\ifnum#1=10 A\else\ifnum#1=11 B\else\ifnum#1=12 C\else
\ifnum#1=13 D\else\ifnum#1=14 E\else\ifnum#1=15 F\fi\fi\fi\fi\fi\fi\fi}
\def\msa@{\hexnumber@\msafam}
\def\llcorner{\delimiter"4\msa@78\msa@78 }
\def\lrcorner{\delimiter"5\msa@79\msa@79 }
\mathchardef\blacktriangleright="3\msa@49
\mathchardef\blacktriangleleft="3\msa@4A
\font\tenmsb=msbm10 scaled 1\@ptsize00
\newfam\msbfam
\textfont\msbfam=\tenmsb
\scriptfont\msbfam=\tenmsb
\def\msb@{\hexnumber@\msbfam}
\mathchardef\varkappa="0\msb@7B
\newdimen\linethick  \linethick=0.4pt
\newdimen\hboxitspace    \hboxitspace=5pt
\newdimen\vboxitspace    \vboxitspace=5pt
\def\fr#1{%
\be\new
\vcenter{
\hrule height\linethick
        \hbox{\vrule width\linethick
              \kern\hboxitspace
              \vbox{\kern\vboxitspace
                    \hbox{$\begin{array}{c}\displaystyle#1
       \end{array}$}%
                    \kern\vboxitspace}%
              \kern\hboxitspace
              \vrule width\linethick}%
        \hrule height\linethick}%
\ee}
\newdimen\Squaresize \Squaresize=14pt
\newdimen\Thickness \Thickness=0.5pt
\def\Square#1{\hbox{\vrule width \Thickness
\vbox to \Squaresize{\hrule height \Thickness\vss
   \hbox to \Squaresize{\hss#1\hss}
\vss\hrule height\Thickness}
\unskip\vrule width \Thickness}
\kern-\Thickness}
\def\Vsquare#1{\vbox{\Square{$#1$}}\kern-\Thickness}

\def\numberbysection{\@addtoreset{equation}{section}
     \def\theequation{\thesection.\arabic{equation}}}
\numberbysection

\renewcommand{\theequation}{\thesection.\arabic{equation}}
\def\titlepage{\@restonecolfalse\if@twocolumn\@restonecoltrue\onecolumn
  \else \newpage \fi \thispagestyle{empty}\c@page\z@
     \def\thefootnote{\fnsymbol{footnote}} }
\def\endtitlepage{\if@restonecol\twocolumn \else  \fi
     \def\thefootnote{\arabic{footnote}}
     \setcounter{footnote}{0}}  
\relax
\hybrid
\makeatletter
\newdimen\normalarrayskip            
\newdimen\minarrayskip               
\normalarrayskip\baselineskip
\minarrayskip\jot
\newif\ifold             \oldtrue            \def\new{\oldfalse}
\def\arraymode{\ifold\relax\else\displaystyle\fi}
\def\eqnumphantom{\phantom{(\theequation)}} 
\def\@arrayskip{\ifold\baselineskip\z@\lineskip\z@
  \else
  \baselineskip\minarrayskip\lineskip1\baselineskip\fi}
\def\@arrayclassz{\ifcase \@lastchclass \@acolampacol \or
\@ampacol \or \or \or \@addamp \or
\@acolampacol \or \@firstampfalse \@acol \fi
\edef\@preamble{\@preamble
\ifcase \@chnum
  \hfil$\relax\arraymode\@sharp$\hfil
  \or $\relax\arraymode\@sharp$\hfil
  \or \hfil$\relax\arraymode\@sharp$\fi}}
\def\@array[#1]#2{\setbox\@arstrutbox=\hbox{\vrule
  height\arraystretch \ht\strutbox
  depth\arraystretch \dp\strutbox
width\z@}\@mkpream{#2}\edef\@preamble{\halign \noexpand\@halignto
\bgroup \tabskip\z@ \@arstrut \@preamble \tabskip\z@ \cr}%
\let\@startpbox\@@startpbox \let\@endpbox\@@endpbox
\if #1t\vtop \else \if#1b\vbox \else \vcenter \fi\fi
\bgroup \let\par\relax
\let\@sharp##\let\protect\relax
\@arrayskip\@preamble}
\def\eqnarray{\stepcounter{equation}%
           \let\@currentlabel=\theequation
           \global\@eqnswtrue
           \global\@eqcnt\z@
           \tabskip\@centering              
           \let\\=\@eqncr
           $$%
         \halign to \displaywidth  \bgroup
          \eqnumphantom \@eqnsel
   \hskip\@centering                               
 $\displaystyle  \tabskip\z@ {##}$%
 &\global\@eqcnt\@ne \hskip 2\arraycolsep
      $ \displaystyle  \arraymode{##}$\hfil
 &\global\@eqcnt\tw@ \hskip 2\arraycolsep
      $\displaystyle\tabskip\z@{##}$\hfil
      \tabskip\@centering
 &{##}\tabskip\z@\cr}
\makeatother

\newtheorem{rem}{Remark}[section]
\newcommand{\beq}[1]{\begin{equation}\label{#1}}
\newcommand\eeq{\end{equation}}
\newcommand\bqa{\begin{eqnarray}}
\newcommand\eqa{\end{eqnarray}}
\def\be{\begin{eqnarray}\new\begin{array}{cc}}
\def\ee{\end{array}\end{eqnarray}}

\def\beq{\begin{equation}}
\def\eeq{\end{equation}}
\def\bse{\begin{subequations}}                
\def\ese{\end{subequations}}
\def\bp{\begin{pmatrix}}
\def\ep{\end{pmatrix}}
\def\bs{\begin{subequations}\label}
\def\es{\end{subequations}}
\def\beq{\begin{eqnarray}{cc}}
\def\eeq{\end{eqnarray}}
\newcommand{\ba}{\begin{align}}
\newcommand{\ea}{\end{align}}
\def\bp{\begin{pmatrix}}
\def\ep{\end{pmatrix}}
\def\bv{\begin{vmatrix}}
\def\ev{\end{vmatrix}}
\def\bel{\be\label}

\def\stack#1#2{\raise0.7pt\hbox{$\mathrel{\mathop{#2}\limits^{#1}}$}}
\def\tr{\triangleright}
\def\tl{\triangleleft}
\def\sem{\mathsurround=0pt \raise1pt
\hbox{$\scriptscriptstyle>\!\!$}\:\!\!\tl}
\def\mes{\mathsurround=0pt \tr\!\:\!\raise0.8pt
\hbox{$\scriptscriptstyle\!\!<$}\,}
\def\]{\mathsurround=0pt ]\raise-2pt\hbox{$_\ast$}}

\def\al{\alpha}

\def\<{\langle}
\def\>{\rangle}

\def\we{\raise-1pt\hbox{$\,\stackrel{\wedge}{,}\,$}}

\begin{document}
\title{
\hfill{\normalsize ITEP-TH-14/14}\\ [10mm] \bf  Theta vocabulary I}
\author{S. Kharchev\thanks{Institute for Theoretical and Experimental
Physics, Moscow, Russia}
\and A. Zabrodin\thanks{Laboratory of 
Mathematical Physics, National Research University
``Higher School of Economics'',
20 Myasnitskaya Ulitsa, Moscow 101000, Russia and
Institute for Theoretical and Experimental
Physics, Moscow, Russia}}
\date{\phantom{.}}
\maketitle

\begin{abstract}
This paper is an annotated list of transformation properties and identities
satisfied by the four theta functions $\theta _1$, $\theta _2$, $\theta _3$, $\theta _4$
of one complex variable, presented in a
ready-to-use form.
An attempt is made to reveal a pattern behind various identities
for the theta-functions.
It is shown that all possible 3, 4 and 5-term identities
of degree four emerge as algebraic
consequences of the six fundamental bilinear 3-term identities
connecting the theta-functions with modular parameters $\tau$ and $2\tau$.
\end{abstract}

\clearpage \newpage
\footnotesize
\normalsize
\newpage
\section{Foreword}

The theta functions introduced by Jacobi \cite{J}
(see also \cite{B}, \cite{W}, \cite{WW}, \cite{M1})
are doubly (quasi)periodic analogues of
the basic trigonometric functions $\sin (\pi u)$ and $\cos (\pi u)$.
Let the two (quasi)periods be $1$ and $\tau \in \CC$ with the condition
$\Im \, \tau >0$.
The basic theta functions are
$\theta _1(u|\tau )$,
$\theta _2(u|\tau )$, $\theta _3(u|\tau )$, $\theta _4(u|\tau )$.
The theory of theta functions is a sort of
``elliptically deformed'' trigonometry.
In essence the functions  $\sin$ and $\cos$ are the same because
$\cos x =\sin (x+\frac{\pi}{2})$, but everybody knows that
in practice it is more convenient to work with
the two functions rather than one.
Likewise, the four theta functions can be obtained from
any one of them by simple transformations like shifts
of the argument and multiplying by a common factor, but it is more
convenient to deal with the set of four instead of one.

\noindent
The ``elliptic deformation'' of the trigonometric functions
may go in two ways depending on which property of the former
one wants to preserve or generalize.
One is a deformation in the class
of {\it entire functions} (the north-east arrow in the diagram below).
It leads to the quasi-periodic theta functions, which are
regular functions in the whole complex plane.
The other one is in the class of {\it doubly periodic functions}.
The (infinite) second period of the
trigonometric functions becomes finite (equal to $\tau$)
at the price of breaking the global analyticity,
so the elliptic functions
$\rm{sn}$, $\rm{cn}$ and $\rm{dn}$, which are
doubly periodic
analogues of trigonometric $\sin$ and $\cos$
are {\it meromorphic functions}
in the complex plane.
$$
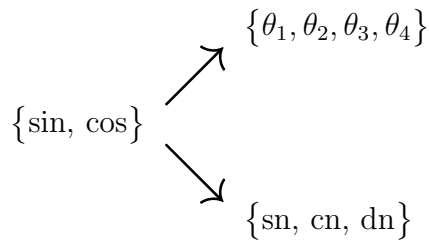
\begin{array}{lll}
& & \!\! \bigl \{ \theta_1, \theta_2, \theta_3, \theta_4 \bigr \}\\
&  \mbox{{\huge $\nearrow$}} & \\
\,\, \bigl \{ \sin , \, \cos \bigr \}\!\! &&\\
& \mbox{{\huge $\searrow$}} &\\
&& \!\! \bigl \{ \rm{sn}, \, {\rm cn }, \, {\rm dn }\bigr \}
\end{array}
$$
In fact the basic elliptic functions
are constructed as ratios of the theta functions and in this sense
the latter seem to be more fundamental.

\noindent
In practical calculations with trigonometric
functions (and their hyperbolic cousins),
one needs just a few identities for the basic functions
$\sin$ and $\cos$ like the addition formula
$\sin (x+y)=\sin x  \cos y + \sin y  \cos x$.
It is not difficult to remember them all
or derive any forgotten one from scratch using the
definitions $\sin x =-i(e^{ix}-e^{-ix})/2$,
$\cos x =(e^{ix}+e^{-ix})/2$.
For the theta functions, the situation is much more involved.
They are connected by a plethora
of identities most of which are not obvious,
not suitable for memorizing and can not be
derived from scratch in any easy way.
Here is what Mumford wrote in Chapter 1
of his book ``Tata lectures on Theta I''
\cite{M1} after presenting a list of ponderous
identities for theta functions:
\begin{itemize}
\item[]
{\it ``We have listed these at such length to illustrate
a key point in the theory
of theta functions: the symmetry of the situation generates rapidly an overwhelming number
of formulae, which do not however make a completely
elementary pattern. To obtain a clear
picture of the algebraic implications of these formulae altogether
is then not usually easy.''}
\end{itemize}

\noindent
All this is aggravated by the fact that
there are several different systems of notation for theta functions
in use.

\noindent
In the present paper we make an attempt to bring some order
into this conglomeration of formulae.
We show that the
3, 4 and 5-term identities of degree four (i.e. with products of
four theta functions in each term), referred to as Weierstrass
addition formulae, Jacobi relations, and Riemann identities,
respectively, can be obtained
by purely algebraic manipulations from {\it six basic
3-term theta relations of degree two}
connecting theta functions with modular parameters $\tau$ and $2\tau$.
Starting with
the six ``elementary bricks'', it is possible to derive
52 fundamental relations of degree four containing four independent
variables.
Besides, we give the complete list of all important particular identities
which are appropriate specifications of the basic
bilinear and degree four ones.
Recently, Koornwinder has proved \cite{K} that
the Weierstrass addition formulae and the Riemann identities
are equivalent. We reproduce this result in a very simple way.

\noindent
In the future we plan to address more specific
questions related to the role of theta functions in the theory
of integrable systems and lattice models of statistical mechanics.

\paragraph{Acknowledgments.}
We are grateful to I. Marshall and
A. Morozov for reading the ma\-nu\-script and
valuable advices.
This work was supported in part by grant NSh-1500.2014.2
for support of scientific schools.
The work of S.K. was supported in part by RFBR grant 15-01-99504.
The work of A.Z. was supported in part by RFBR grants
15-01-05990 and 14-01-90405-Ukr.
The financial support from the Government of the Russian Federation 
within the framework of the implementation of the 5-100 Programme Roadmap of 
the National Research University  Higher School of Economics is acknowledged.

\section{Theta functions}\label{sec1}
\subsection{Theta functions with characteristics}
Fix the modular parameter $\tau \in \CC$ such that $\Im \, \tau >0$
and consider the infinite series \cite{M1}:
\bel{ch1}
\theta_{a,b}(u|\tau)=\sum_{k\in\ZZ}
\exp\Big\{\pi i\tau(k+a)^2+2\pi i(k+a)(u+b)\Big\},
\ee
where $i=\sqrt{-1}$ and $a,b\in\RR$. The series is absolutely convergent for
any $u\in \CC$ and defines
the entire function $\theta_{a,b}(u|\tau)$.
It is called the theta function with characteristics $a,b$.
These functions
are connected by the relations
\bs{ch3}
\ba
\theta_{a,b}(u+\tau a'+b'|\tau )&=e^{-2\pi ia'(u+b+b'+a'\tau/2)}
\theta_{a+a',b+b'}(u|\tau),\label{ch3a}\\
\theta_{a+1,b}(u|\tau)&=\theta_{a,b}(u|\tau)\label{ch3b},\\
\theta_{a,b+1}(u|\tau)&=e^{2\pi i a}\theta_{a,b}(u|\tau ).\label{ch3c}
\end{align}
\es
In particular, it follows from here that
the functions $\theta_{a,b}$ are quasiperiodic with
(quasi)periods $1$ and $\tau$:
\bs{ch7}
\ba
\label{ch7a}
\theta_{a,b}(u+1|\tau )&=e^{2\pi i a}\theta_{a,b}(u|\tau ),\\
\label{ch7b}
\theta_{a,b}(u+\tau |\tau )&=e^{-\pi i(2u+2b+\tau)}\theta_{a,b}(u|\tau ),
\end{align}
\es
and the shifts by half-periods are given by
\bs{ch6}
\ba
&\theta_{a,b}(u+{\textstyle\frac{1}{2}}|\tau)
=\theta_{a,b+\frac{1}{2}}(u|\tau )=
e^{2\pi i a}\theta_{a,b-\frac{1}{2}}(u|\tau ),\label{ch6a}\\
&\theta_{a,b}(u+{\textstyle\frac{\tau}{2}}|\tau)=
e^{-\pi i(u+b+\tau/4)}
\theta_{a+\frac{1}{2},b}(u|\tau)=
e^{-\pi i(u+b+\tau/4)}\theta_{a-\frac{1}{2},b}(u|\tau).
\label{ch6b}
\end{align}
\es
It follows from the definition that $\theta_{a,b}(-u)=\theta_{-a,-b}(u)$.
Hence, according to (\ref{ch3b}), (\ref{ch3c}),
the functions $\theta_{a,b}(u)$ have definite evenness properties only for integer or half-integer characteristics.
By virtue of (\ref{ch3b}), (\ref{ch3c}), it is sufficient
to consider the theta functions with characteristics $0\leq a,b<1$.

\subsection{Basic theta functions}
The standard theta functions with half-integer characteristics \cite{J,B}
are defined as follows:
\bel{d1}
\begin{array}{l}
\theta_1(u|\tau )=-\,\theta_{\frac{1}{2},\frac{1}{2}}(u|\tau)
=-i\sum_{k\in\ZZ}(-1)^kq^{(k+\frac{1}{2})^2}
e^{\pi i(2k+1)u},\\
\theta_2(u|\tau )=\,\theta_{\frac{1}{2},0}(u|\tau)
=\,\, \sum_{k\in\ZZ}q^{(k+\frac{1}{2})^2}e^{\pi i(2k+1)u},\\
\theta_3(u|\tau )= \, \theta_{0,0}(u|\tau)
=\,\, \sum_{k\in\ZZ}q^{k^2}e^{2\pi iku},\\
\theta_4(u|\tau )=\, \theta_{0,\frac{1}{2}}(u|\tau)
=\,\, \sum_{k\in\ZZ}(-1)^kq^{k^2}e^{2\pi iku},
\end{array}
\ee
where
\bel{d3}
q:=e^{\pi i\tau},\ \ \ |q|<1.
\ee
In the limit $\tau \to i\infty$ they are:
$\theta_1(u|\tau )=2q^{\frac{1}{4}}\sin \pi u + O(q^{\frac{9}{4}})$,
$\theta_2(u|\tau )=2q^{\frac{1}{4}}\cos \pi u + O(q^{\frac{9}{4}})$,
$\theta_3(u|\tau )=1+O(q)$, $\theta_4(u|\tau )=1+O(q)$.

\noindent
In what follows we often write $\theta_r(u|\tau ):=\theta_r(u),\,r=1,2,3,4$
if this does not cause confusion. From (\ref{d1}) it is clear that the function
$\theta_1$ is odd, $\theta _1(-u)=-\theta_1(u)$; the other three are even,
$\theta_s(-u)=\theta_s(u)$, $s=2,3,4$.

\noindent
The values $\theta_1'(0)$, $\theta_2(0)$, $\theta_3(0)$, $\theta_4(0)$ are called
theta constants.

\subsection{Shifts by periods and half-periods}
Here we list the essential transformation properties for the theta functions
(\ref{d1}) which follow from (\ref{ch3}).

\smallskip

\noindent
\underline{Shifts by periods}:
\bel{p}
\begin{array}{l}
\theta_1(u+1)=-\theta_1(u),\\
\theta_2(u+1)=-\theta_2(u),\\
\theta_3(u+1)=\theta_3(u),\\
\theta_4(u+1)=\theta_4(u).
\end{array}\hspace{2cm}
\begin{array}{l}
\theta_1(u+\tau)=-e^{-\pi i(2u+\tau)}\theta_1(u),\\
\theta_2(u+\tau)=e^{-\pi i(2u+\tau)}\theta_2(u),\\
\theta_3(u+\tau)=e^{-\pi i(2u+\tau)}\theta_3(u),\\
\theta_4(u+\tau)=-e^{-\pi i(2u+\tau)}\theta_4(u).
\end{array}\vspace{0.3cm}\\
\begin{array}{l}
\theta_1(u \! +\! \tau \! +\! 1)=e^{-\pi i(2u+\tau)}\theta_1(u),\\
\theta_2(u\! +\! \tau \! +\! 1)=-e^{-\pi i(2u+\tau)}\theta_2(u),\\
\theta_3(u\! +\! \tau \! +\! 1)=e^{-\pi i(2u+\tau)}\theta_3(u),\\
\theta_4(u\! +\! \tau \! +\! 1)=-e^{-\pi i(2u+\tau)}\theta_4(u).
\end{array}
\ee
\underline{Shifts by half-periods}:
\bel{hp}
\begin{array}{l}
\theta_1(u+{\textstyle\frac{1}{2}})=\theta_2(u),\\
\theta_2(u+{\textstyle\frac{1}{2}})=-\theta_1(u),\\
\theta_3(u+{\textstyle\frac{1}{2}})=\theta_4(u),\\
\theta_4(u+{\textstyle\frac{1}{2}})=\theta_3(u).
\end{array}\hspace{2cm}
\begin{array}{l}
\theta_1(u+{\textstyle\frac{\tau}{2}})=
ie^{-\pi i(u+\tau/4)}\theta_4(u),\\
\theta_2(u+{\textstyle\frac{\tau}{2}})=
e^{-\pi i(u+\tau/4)}\theta_3(u),\\
\theta_3(u+{\textstyle\frac{\tau}{2}})=
e^{-\pi i(u+\tau/4)}\theta_2(u),\\
\theta_4(u+{\textstyle\frac{\tau}{2}})=i
e^{-\pi i(u+\tau/4)}\theta_1(u).
\end{array}\vspace{0.3cm}\\
\begin{array}{l}
\theta_1\bigl (u+{\textstyle\frac{\tau+1}{2}}\bigr )=
e^{-\pi i(u+\tau/4)}\theta_3(u),\\
\theta_2 \bigl (u+{\textstyle\frac{\tau+1}{2}}\bigr )=-i e^{-\pi i(u+\tau/4)}\theta_4(u),\\
\theta_3 \bigl (u+{\textstyle\frac{\tau+1}{2}}\bigr )=
ie^{-\pi i(u+\tau/4)}\theta_1(u),\\
\theta_4\bigl (u+{\textstyle\frac{\tau+1}{2}}\bigr )=
e^{-\pi i(u+\tau/4)}\theta_2(u).
\end{array}
\ee

\subsection{Zeros of theta functions}
These relations imply that the (first order) zeros of the theta
functions are as follows:
\begin{equation}\label{zeros}
\begin{array}{l}
\theta_1(u)=0: \quad u=n+m\tau,
\\
\theta_2(u)=0: \quad u=n+\frac{1}{2}+m\tau,
\\
\theta_3(u)=0: \quad u=n+\frac{1}{2} +(m+\frac{1}{2})\tau,
\\
\theta_4(u)=0: \quad u=n +(m+\frac{1}{2})\tau,
\end{array}
\end{equation}
where $n,m\in \ZZ$.

\smallskip

{\small
\noindent
Indeed, in accordance with (\ref{p}), (\ref{hp}),
the functions $\theta_1(u  + n + m\tau)$,
$\theta_2 \bigl (u + n + \frac{1}{2} + m\tau \bigr )$,
$\theta_3 \bigl (u\! +\! n\! +\! \frac{1}{2} \! +\!
(m \! +\! \frac{1}{2})\tau \bigr )$, $\theta_1 \bigl (u\! +\! n\!  +\!
(m \! +\! \frac{1}{2})\tau \bigr )$
are proportional to the odd function $\theta_1(u)$.
Hence the corresponding zeros are as in (\ref{zeros}).
To complete the proof,
it is sufficient to show that the function
$\theta_1(u)$ has precisely one simple zero
in the parallelogram $\Pi$ with the vertices
$-\frac{1}{2}\pm \frac{\tau}{2}, \frac{1}{2}\pm \frac{\tau}{2}$. The
standard argument is to compute the contour integral
$\frac{1}{2\pi  i}\oint_{\partial \Pi} d\log\theta_1(u)=1$ which means that
the zero is simple
(see \cite{WW}, \cite{M1} for details).
}

\subsection{Theta functions as infinite products}
One has the following infinite product representations:
\bs{d4}
\ba
\theta_1(u|\tau)&=2q^{\frac{1}{4}}\sin\pi u
\prod_{n=1}^\infty(1-q^{2n})(1-q^{2n}e^{2\pi i u})(1-q^{2n}e^{-2\pi i u}),\label{d4a}\\
\theta_2(u|\tau)&=2q^{\frac{1}{4}}\cos\pi u
\prod_{n=1}^\infty(1-q^{2n})(1+q^{2n}e^{2\pi i u})(1+q^{2n}
e^{-2\pi i u}),\label{d4b}\\
\theta_3(u|\tau)&=
\prod_{n=1}^\infty(1-q^{2n})(1+q^{2n-1}e^{2\pi i u})(1+q^{2n-1}e^{-2\pi i u}),
\label{d4c}\\
\theta_4(u|\tau)&=
\prod_{n=1}^\infty(1-q^{2n})(1-q^{2n-1}e^{2\pi i u})(1-q^{2n-1}e^{-2\pi i u})
\label{d4d}.
\end{align}
\es
{\small
To prove (\ref{d4d}), we note that the product
$p(u|\tau):=\prod_{n=1}^\infty(1-q^{2n-1}e^{2\pi i u})
(1-q^{2n-1}e^{-2\pi i u})$
has the same zeros as $\theta_4(u|\tau)$ and the ratio $\theta_4(u|\tau)/p(u|\tau)$
is a doubly periodic function with periods 1 and $\tau$. Hence the ratio is constant and one has $\theta_4(u|\tau)=Ap(u|\tau)$.
To find the constant $A$,
put $u=0$ thus
getting $A=\theta_4(0|\tau)/p(0|\tau)$. Finally,
in accordance with the Gauss formula
\cite[p. 23, eq. (2.2.12)]{An},
\bel{g1}
\theta_4(0|\tau)=\sum_{k\in\ZZ}(-1)^kq^{k^2}=
\prod_{n=1}^\infty\frac{1-q^n}{1+q^n}.
\ee
Rewriting $\prod_{n=1}^\infty(1-q^{2n-1})=
\prod_{n=1}^\infty(1-q^n)/(1-q^{2n})$, one gets
$A=\prod_{n=1}^\infty(1-q^{2n})$ and
formula (\ref{d4d}) is proved.
Equations (\ref{d4a})--(\ref{d4c})
are obtained from (\ref{d4c}) by appropriate shifts
of $u$ in accordance with (\ref{hp}).
}

\noindent
As a corollary of (\ref{d4}) one has infinite product
representations for the
theta constants:
\bs{dr1}
\ba
\theta'_1(0)&=2\pi q^{\frac{1}{4}}\prod_{n=1}^\infty(1-q^{2n})^3,\\
\theta_2(0)&=2q^{\frac{1}{4}}
\prod_{n=1}^\infty(1-q^{2n})(1+q^{2n})^2,\\
\theta_3(0)&=\prod_{n=1}^\infty(1-q^{2n})(1+q^{2n-1})^2,\\
\theta_4(0)&=\prod_{n=1}^\infty(1-q^{2n})(1-q^{2n-1})^2.
\end{align}
\es
Since
$\displaystyle{\prod_{n\geq 1}(1+q^{2n})(1+q^{2n-1})(1-q^{2n-1})=1}$,
this implies the famous identity for the theta constants \cite[p. 517]{J}:
\bel{tc1}
\theta_1'(0)=\pi\theta_2(0)\theta_3(0)\theta_4(0).
\ee

\subsection{Modular transformations}

\noindent
\underline{The transformation $\tau \to \tau +1$}:
\bs{mod1}
\ba
\theta_1 (u|\tau+1)&=e^{\frac{\pi i}{4}}\theta_1 (u|\tau ),\label{mod1a}\\
\theta_2 (u|\tau+1)&=e^{\frac{\pi i}{4}}\theta_2 (u|\tau),\label{mod1b}\\
\theta_3 (u|\tau+1)&=\theta_4(u|\tau),\label{mod1c}\\
\theta_4 (u|\tau+1)&=\theta_3(u|\tau).\label{mod1d}
\end{align}
\es
Since $\tau\rightarrow\tau+1$ implies $q\rightarrow -q$, equations
(\ref{mod1}) follow from (\ref{d1}) or (\ref{d4}).
\smallskip

\noindent
\underline{The transformation $\tau \to -1/\tau$}:
\bs{mod2}
\ba
\theta_1\left (u/\tau|-1/\tau\right)&=-i\sqrt{-i\tau}\,
e^{\pi iu^2/\tau}\theta_1(u|\tau),\label{mod2a}\\
\theta_2\left (u/\tau|-1/\tau\right)&=\sqrt{-i\tau}\,
e^{\pi iu^2/\tau}\theta_4(u|\tau),\label{mod2b}\\
\theta_3\left (u/\tau|-1/\tau\right)&=\sqrt{-i\tau}\,
e^{\pi iu^2/\tau}\theta_3(u|\tau),\label{mod2c}\\
\theta_4\left (u/\tau|-1/\tau\right)&=\sqrt{-i\tau}\,
e^{\pi iu^2/\tau}\theta_2(u|\tau).\label{mod2d}
\end{align}
\es
The branch of the square root here
is such that $\Re \sqrt{-i\tau}>0$.

\noindent
{\small The proof is well known. Since the ratio
$e^{-\pi iu^2/\tau}\theta_3
\left(u/\tau |-1/\tau \right)/\theta_3(u|\tau):=C$ is an entire
doubly periodic function of $u$, it is a constant (which may depend
only on $\tau$). Shifting $u$ by
$\frac{1}{2}, \frac{\tau}{2}$, $\frac{\tau+1}{2}$, one obtains
three more relations of the same kind with
the same constant $C$. Then the substitution of these formulas
to (\ref{tc1}), yields $C^2=-i\tau$.
The sign of the square root is determined
by the argument that if $\tau \in i\RR_{+}$, then both
$\theta_3(0|-1/\tau )$ and $\theta_3(0|\tau)$ are real and positive.
}

\subsection{Other notation for the theta functions}
The notations for theta functions used in the literature
are of a great variety.
This can be a source of confusion. Here we briefly comment on the main
systems of notation other that the one adopted in this paper.
In the theory of
elliptic integrals,
the theta functions
\bel{Thetas}
\Theta _r \bigl (u|\tau \bigr )=\theta_r\Bigl(
\frac{u}{2K}\Bigm | \tau \Bigr ), \quad
K=\frac{\pi}{2}\, \theta_3^2(0|\tau )\,
\ee
introduced by Riemann are commonly encountered.
The number $K$ is the full elliptic integral (of the
first kind).
In \cite{A} and some other books
our $\theta_r$ is denoted as $\vartheta_r$ while
$\Theta_r$ defined in (\ref{Thetas}) is just $\theta_r$.
The antiquated Jacobi notation (still preferred by some authors) are
$H$, $H_1$, $\Theta_1$, $\Theta$ for $\Theta_1$,
$\Theta_2$, $\Theta_3$, $\Theta_4$ respectively.
The ``multiplicative notation'' $\theta_r(z|q)$ for $\theta_r(u|\tau )$,
where $q=e^{\pi i\tau}$, $z=e^{2\pi i u}$, is
widely used in the modern literature on elliptic hypergeometric
series and related problems.

\noindent
Lastly, let us mention
a few of the minor differences in notation encountered in the literature.
In {\rm\cite{W}}, {\rm\cite{HC}} the functions $\Theta_{a,b}(u)$ have been
considered which are related with $\theta_{a,b}(u)$ by
\bel{ch4}
\Theta^{\rm\scriptscriptstyle W}_{a,b}(u)=
e^{\pi iab}\theta_{-\frac{a}{2},\frac{b}{2}}(u), \quad
\Theta^{\rm\scriptscriptstyle HC}_{a,b}(u)=
e^{-\frac{\pi iab}{2}}\theta_{\frac{a}{2},\frac{b}{2}}(u).
\ee
The set of our theta functions (\ref{d1}) is related with the corresponding
functions in {\rm \cite{WW}} as
$\theta_r(u|\tau)=\theta^{{\rm\scriptscriptstyle WW}}_r(\pi u|\tau),\;r=1, 2, 3, 4$.
Following the original notation \cite{J}, in \cite{A}
and in some other books the notation
$\theta_0$ is used instead of $\theta_4$.

\section{Four types of identities between theta functions}
\subsection{Preliminaries}

The number of identities satisfied by the theta functions is enormous.
It is still fairly big if we consider
identities involving up to four independent variables.
They can be split into four types:
\begin{itemize}
\item[]
{\bf B.} Three-term bilinear identities involving two independent
variables.
They relate products of
two theta functions with modular parameter $\tau$ to linear
combinations (actually, sums or differences) of similar
products of theta functions with modular parameter $2\tau$.\\
{\bf W.}
Three-term identities of degree 4 (the
{\it Weierstrass addition formulae}).\\
{\bf J.} Four-term identities  of degree 4 (the {\it Jacobi formulae}).\\
{\bf R.} Five-term identities of degree 4 (the {\it Riemann identities}).
\end{itemize}

\noindent
The identities of types {\bf W}, {\bf J}, {\bf R} include theta functions
with the same modular parameter $\tau$ and contain four independent
variables.
The identities of type {\bf B} are the most fundamental ones:
all the others are algebraic consequences of these
together with the evenness properties of the theta functions
$\theta_r(-u)=(-1)^{\delta_{r,1}}\theta_r(u),\, r=1,2,3,4$.
Namely, we shall show how to derive {\bf W} from {\bf B} etc.,
according to the
scheme {\bf B}\,$\to$\,{\bf W}\,$\to$\,{\bf J}\,$\to$\,{\bf R}.
It also turns out
that the Jacobi and Riemann identities are equivalent in a
very simple way.
At the end of this
section, we prove the arrow ${\bf W}\leftarrow{\bf J}$ which
implies equivalence of the Weierstrass and Jacobi identities.

\subsection{Three-term bilinear identities connecting theta functions with
$\tau$ and $2\tau$}\label{bil}
{\bf B.I}. There are six basis bilinear identities:
\bs{sb1}
\ba
\theta_1(u|\tau)\theta_1(v|\tau)=\theta_3(u+v|2\tau)\theta_2(u-v|2\tau)-
\theta_2(u+v|2\tau)\theta_3(u-v|2\tau)\label{sb1a},\\
\theta_1(u|\tau)\theta_2(v|\tau)=\theta_1(u+v|2\tau)\theta_4(u-v|2\tau)+
\theta_4(u+v|2\tau)\theta_1(u-v|2\tau)\label{sb1b},\\
\theta_2(u|\tau)\theta_2(v|\tau)=\theta_2(u+v|2\tau)\theta_3(u-v|2\tau)+
\theta_3(u+v|2\tau)\theta_2(u-v|2\tau)\label{sb1c},\\
\theta_3(u|\tau)\theta_3(v|\tau)=\theta_3(u+v|2\tau)\theta_3(u-v|2\tau)+
\theta_2(u+v|2\tau)\theta_2(u-v|2\tau)\label{sb1d},\\
\theta_3(u|\tau)\theta_4(v|\tau)=\theta_4(u+v|2\tau)\theta_4(u-v|2\tau)-
\theta_1(u+v|2\tau)\theta_1(u-v|2\tau)\label{sb1e},\\
\theta_4(u|\tau)\theta_4(v|\tau)=\theta_3(u+v|2\tau)\theta_3(u-v|2\tau)-
\theta_2(u+v|2\tau)\theta_2(u-v|2\tau)\label{sb1f}.
\end{align}
\es
(See \cite{Ig}, \cite{D}, \cite{M2} for the general case of multi-dimensional
theta functions.)

\noindent
{\bf B.II}. A system equivalent to (\ref{sb1}):
\bs{es1}
\ba
2\theta_1(u+v|2\tau)\theta_1(u-v|2\tau)&=\theta_4(u|\tau)\theta_3(v|\tau)-
\theta_3(u|\tau)\theta_4(v|\tau)\label{es1a},\\
2\theta_1(u+v|2\tau)\theta_4(u-v|2\tau)&=\theta_1(u|\tau)\theta_2(v|\tau)+
\theta_2(u|\tau)\theta_1(v|\tau)\label{es1b},\\
2\theta_2(u+v|2\tau)\theta_2(u-v|2\tau)&=\theta_3(u|\tau)\theta_3(v|\tau)-
\theta_4(u|\tau)\theta_4(v|\tau)\label{es1c},\\
2\theta_2(u+v|2\tau)\theta_3(u-v|2\tau)&=\theta_2(u|\tau)\theta_2(v|\tau)-
\theta_1(u|\tau)\theta_1(v|\tau)\label{es1d},\\
2\theta_3(u+v|2\tau)\theta_3(u-v|2\tau)&=\theta_3(u|\tau)\theta_3(v|\tau)+
\theta_4(u|\tau)\theta_4(v|\tau)\label{es1e},\\
2\theta_4(u+v|2\tau)\theta_4(u-v|2\tau)&=\theta_3(u|\tau)\theta_4(v|\tau)+
\theta_4(u|\tau)\theta_3(v|\tau)\label{es1f}.
\end{align}
\es
\begin{rem}
Starting with any identity in (\ref{sb1}),
one can derive all the other ones
by appropriate shifts of the variables $u$, $v$.
\end{rem}
{\small
The proof is standard. Let us
prove, for example, (\ref{sb1b}). Consider the function
$$
F(v):=\theta_1(u+v|2\tau)\theta_4(u-v|2\tau)+\theta_4(u+v|2\tau)\theta_1(u-v|2\tau).
$$
By virtue of (\ref{p}) and (\ref{hp}),
$F(v+1)=-F(v)$, $F(v+\tau)=e^{-\pi i(2v+\tau)}F(v)$
and $F(\frac{1}{2})=0$. Hence
zeros of $F(v)$ are $v_{n,m}=n+\frac{1}{2}+m\tau$, $n,m\in\ZZ$ and the ratio
$F(v)/(\theta_1(u|\tau)\theta_2(v|\tau))$ is an entire
function doubly periodic in $v$ with periods
$1, \tau$.
Therefore, this ratio does not depend on $v$:
$F(v)/\theta_1(u|\tau)\theta_2(v|\tau)=C(u)$. Setting $v=u$, one has:
$$
C(u)=\frac{\theta_1(2u|2\tau)\, \theta_4(0|2\tau)}{\theta_1(u|\tau)\,
\theta_2(u|\tau)}\,.
$$
By virtue of (\ref{d4a}), (\ref{d4b}) and (\ref{g1}),
$\theta_1(u|\tau)\theta_2(u|\tau)=\theta_1(2u|2\tau)\theta_4(0|2\tau)$
and thus $C(u)\equiv 1$.
}

\subsection{Three-term Weierstrass addition identities}\label{qr3}
There are twelve addition formulae (see below).
We start with the identity
\bel{tt0}
\hspace{-1cm}
\theta_1(u+x)\theta_1(u-x)\theta_1(v+y)\theta_1(v-y)-
\theta_1(u+y)\theta_1(u-y)\theta_1(v+x)\theta_1(v-x)\\=
\theta_1(u+v)\theta_1(u-v)\theta_1(x+y)\theta_1(x-y)
\ee
which was originally discovered and
proved by Weierstrass \cite[p. 155]{We}.
All the
identities listed below in this section
can be derived from it by appropriate shifts
of the variables in accordance with relations (\ref{hp}).
Our approach is different. We show that
all the identities of Weierstrass' type are simple algebraic consequences of
the bilinear system (\ref{sb1}) together
with the evenness conditions
$\theta_r(-u)=(-1)^{\delta_{r,1}}\theta_r(u)$, $r=1,2,3,4$.
This argument is independent of (\ref{hp}).

{\small
\smallskip
\noindent
To prove (\ref{tt0}), one should rewrite (\ref{sb1a}) as
$$
\theta_1(u+x|\tau)\theta_1(u-x|\tau)=
\theta_3(2u|2\tau)\theta_2(2x|2\tau)-
\theta_2(2u|2\tau)\theta_3(2x|2\tau).
$$
Multiply this by the similar
expression for
$\theta_1(v+y|\tau)\theta_1(v-y|\tau)$
and subtract the same with the change
$x\leftrightarrow y$. Using (\ref{sb1a}) once again,
we arrive at (\ref{tt0}). All the equations
below in this section can be obtained from system (\ref{sb1})
in a similar way.
}

\noindent
{\bf W.I}. Symmetric system:
\bel{tt1}
\hspace{-0.5cm}
\theta_1(u+x)\theta_1(u-x)\theta_r(v+y)\theta_r(v-y)-
\theta_1(v+x)\theta_1(v-x)\theta_r(u+y)\theta_r(u-y)\\=
\theta_1(u+v)\theta_1(u-v)\theta_r(x+y)\theta_r(x-y),
\ee
$r=1,2,3,4$.
\smallskip

\noindent
{\bf W.II}. Complimentary system:
\bse\label{tt2}
\bel{tt2c}
\hspace{-1.3cm}
\theta_2(u+x)\theta_2(u-x)\theta_3(v+y)\theta_3(v-y)-
\theta_2(v+x)\theta_2(v-x)\theta_3(u+y)\theta_3(u-y)\\=-
\theta_1(u+v)\theta_1(u-v)\theta_4(x+y)\theta_4(x-y),
\ee
\bel{tt2b}
\hspace{-1.4cm}
\theta_2(u+x)\theta_2(u-x)\theta_4(v+y)\theta_4(v-y)-
\theta_2(v+x)\theta_2(v-x)\theta_4(u+y)\theta_4(u-y)\\=-
\theta_1(u+v)\theta_1(u-v)\theta_3(x+y)\theta_3(x-y),
\vspace{-0.2cm}
\ee
\bel{tt2a}
\hspace{-1.4cm}
\theta_3(u+x)\theta_3(u-x)\theta_4(v+y)\theta_4(v-y)-
\theta_3(v+x)\theta_3(v-x)\theta_4(u+y)\theta_4(u-y)\\=-
\theta_1(u+v)\theta_1(u-v)\theta_2(x+y)\theta_2(x-y).
\vspace{-0.2cm}
\ee
\ese
\smallskip

\noindent
{\bf W.III}. Asymmetric system:
\bel{tt3}
\hspace{-1cm}
\theta_r(u+x)\theta_r(u-x)\theta_r(v+y)\theta_r(v-y)-
\theta_r(u+y)\theta_r(u-y)\theta_r(v+x)\theta_r(v-x)\\=\,
(-1)^{r-1}\theta_1(u+v)\theta_1(u-v)\theta_1(x+y)\theta_1(x-y),
\ee
$r=1,2,3,4$.
\smallskip

\noindent
{\bf W.IV}. Complimentary identity:
\bel{tt4}
\hspace{-1cm}
\theta_3(u+x)\theta_3(u-x)\theta_3(v+y)\theta_3(v-y)-
\theta_4(v+x)\theta_4(v-x)\theta_4(u+y)\theta_4(u-y)\\=\,
\theta_2(u+v)\theta_2(u-v)\theta_2(x+y)\theta_2(x-y).
\ee
\smallskip

\noindent
{\bf W.V}. Mixed identity:
\bel{tt5}
\hspace{-1cm}
\theta_1(u+x)\theta_2(u-x)\theta_3(v+y)\theta_4(v-y)-
\theta_1(u-y)\theta_2(u+y)\theta_3(v-x)\theta_4(v+x)\\=\,
\theta_1(x+y)\theta_2(x-y)\theta_3(u+v)\theta_4(u-v).
\ee
\begin{rem}
Sometimes the Weierstrass addition formula (\ref{tt0}) is referred to as
Fay identity. In fact, it is a
generalization of Jacobi's results (see Section \ref{paf} below).
\end{rem}

\subsection{Four-term Jacobi identities}\label{qr4}
In Sect. \ref{qr3} we have
presented twelve three-term identities of degree four
depending on four variables $u,\,v,\,x,\,y$.
Here we introduce another set of
variables $W,\,X,\,Y,\,Z$ and their {\it ``dual''} counterparts \cite{WW}:
\vspace{-0.2cm}
\bel{ft1}
\vspace{-0.06cm}
\ \ W'=\textstyle{\frac{1}{2}}(-W+X+Y+Z),\\
\vspace{-0.05cm}
X'=\textstyle{\frac{1}{2}}(W-X+Y+Z),\\
\vspace{-0.06cm}
Y'=\textstyle{\frac{1}{2}}(W+X-Y+Z),\\
Z'=\textstyle{\frac{1}{2}}(W+X+Y-Z).
\ee
One can easily verify that $W,X,Y,Z$ are
expressed via the ``dual'' variables $W',\,X',\,Y',\,Z'$
by the same formulae, i.e.,
the ``prime procedure'' applied to (\ref{ft1}) yields
$(W')'=W=\textstyle{\frac{1}{2}}(-W'+X'+Y'+Z')$ etc.
We employ the short-hand notation
$$[pqrs]:=\theta_p(W)\theta_q(X)\theta_r(Y)\theta_s(Z),
\qquad
[pqrs]':=\theta_p(W')\theta_q(X')\theta_r(Y')\theta_s(Z')
$$
which is widely used in  \cite{WW}.
If all the indices of
the theta functions coincide, this is further abbreviated to
$$[r]:=\theta_r(W)\theta_r(X)\theta_r(Y)\theta_r(Z),
\qquad
[r]':=\theta_r(W')\theta_r(X')\theta_r(Y')\theta_r(Z').
$$
Below we list all four-term basic identities of
degree four which were essentially
obtained by Jacobi \cite[p. 507]{J}.
Here we present these in a more symmetric and
comprehensive form.

\noindent
The simplest (and most important) ones are:
\bs{ft2}
\ba
\hspace{-2.7cm}
[1]+[2] & = [1]'+[2]'\label{ft2a},\\
\hspace{-5cm}
{\bf J.I}\hspace{2cm}
[1]-[2] &= [4]'-[3]'\label{ft2b},\\
[3]+[4] &=[3]'+[4]'\label{ft2c},\\
[3]-[4] &=[2]'-[1]'\label{ft2d}.
\end{align}
\es
{\small
The system (\ref{ft2}) is a direct algebraic
corollary of appropriate addition formulae given in Section \ref{qr3}.
To see this, we relate the variables $\,u,\,v,\,x,\,y$
with the variables of
the present section as follows:
\bel{ft'}
\left\{
\begin{array}{l}
W=u+x,\vspace{-0.2cm}\\
X=u-x,\vspace{-0.2cm}\\
Y=v+y,\vspace{-0.2cm}\\
Z=v-y.
\end{array}
\right.\hspace{0.5cm}\Longleftrightarrow\hspace{0.5cm}
\left\{\begin{array}{l}
W'=v-x,\vspace{-0.2cm}\\
X'=v+x,\vspace{-0.2cm}\\
Y'=u-y,\vspace{-0.2cm}\\
Z'=u+y.
\end{array}\right.
\ee
Further, the products of theta functions
containing ``inappropriate'' combinations $u\pm v,\,x\pm y$
can be excluded from the addition formulae.
Then identities (\ref{ft2a}), (\ref{ft2c}) emerge as particular cases
of (\ref{tt3}).
Changing $v\leftrightarrow x$ in (\ref{tt3}) (with $r=3$) and (\ref{tt4}),
one obtains (\ref{ft2b}). Finally, (\ref{ft2d}) is a ``dual'' version of (\ref{ft2b}).
}
\begin{rem}Equations (\ref{ft2}) are a part of the system of twelve identities
written in {\rm \cite[pp. 468, 488]{WW}}. It is easy to see that all additional
relations are appropriate linear combinations of the basic ones, (\ref{ft2}).
For completeness, we give here the full list:
\bel{ft3}
[1]+[2]=[1]'+[2]', \ \ \ \ [1]+[3]=[2]'+[4]', \ \ \ \ [1]+[4]=[1]'+[4]',\\ \,\!
[2]+[3]=[2]'+[3]',\ \ \ \ [2]+[4]=[1]'+[3]', \ \ \ \ [3]+[4]=[3]'+[4]',\\ \,\!
[1]-[2]=[4]'-[3]', \ \ \ \ [1]-[3]=[1]'-[3]', \ \ \ \ [1]-[4]=[2]'-[3]',\\ \,\!
[2]-[3]=[1]'-[4]',\ \ \ \ [2]-[4]=[2]'-[4]',\ \ \ \ [3]-[4]=[2]'-[1]'.
\ee
\end{rem}

\noindent
Now we list symmetric self-dual identities for products of type $[rrss]$
which can also be
derived from the addition formulae by algebraic manipulations:
\bs{si1}
\begin{align}
[1122]+[2211]=[1122]'+[2211]',\label{si1a}\\ \,
[1133]+[3311]=[1133]'+[3311]',\label{si1b}\\ \,
\hspace{-4cm} {\bf J.II} \hspace{2cm}
[1144]+[4411]=[1144]'+[4411]',\label{si1c}\\ \,
[2233]+[3322]=[2233]'+[3322]',\label{si1d}\\ \,
[2244]+[4422]=[2244]'+[4422]',\label{si1e}\\ \,
[3344]+[4433]=[3344]'+[4433]'\label{si1f}.
\end{align}
\es
Further, there are simple complimentary relations:
\bs{si2}
\begin{align}
[1122]-[2211]=[3344]'-[4433]',\label{si2a}\\ \,
[1133]-[3311]=[2244]'-[4422]',\label{si2b}\\ \,
\hspace{-4cm} {\bf J.III} \hspace{2cm}
[1144]-[4411]=[2233]'-[3322]',\label{si2c}\\ \,
[2233]-[3322]=[1144]'-[4411]',\label{si2d}\\ \,
[2244]-[4422]=[1133]'-[3311]',\label{si2e}\\ \,
[3344]-[4433]=[1122]'-[2211]'\label{si2f}.
\end{align}
\es
The subsystems (\ref{si2a})--(\ref{si2c})
and (\ref{si2d})--(\ref{si2f}) are dual to
each other. Thus, the systems (\ref{si1}), (\ref{si2}) can be represented in very
compact form:
\bel{si3}
\hspace{0.7cm}
[rrss]+[ssrr]=[rrss]'+[ssrr]',\\ \!\,
\hspace{0.7cm}
[rrss]-[ssrr]=[\tilde r\tilde r\tilde s\tilde s]'-
[\tilde s\tilde s\tilde r\tilde r]',
\ee
where $r,s\in(1,2,3,4),\,r<s$ and
$\tilde s,\tilde r\in(1,2,3,4)\backslash(r,s),\ \tilde r<\tilde s$.

\smallskip\noindent
Finally, there are four ``fully mixed'' identities:
\bs{si4}
\ba
[1234]+[2143]=[3412]'+[4321]'\label{ft11a},\\ \,\!
\hspace{-4cm} {\bf J.IV} \hspace{2cm}
[1234]-[2143]=[2143]'-[1234]'\label{ft11b},\\ \,\!
[3412]+[4321]=[1234]'+[2143]'\label{ft11c},\\
[3412]-[4321]=[4321]'-[3412]'\label{ft11d}.
\end{align}
\es

\noindent
{\small
Identities (\ref{si1a})--(\ref{si1c}) follow from (\ref{tt1}).
Indeed, one can write (\ref{tt1}) as
$$
[11rr]-[11rr]'=\theta_1(u+v)\theta_1(u-v)\theta_r(x+y)\theta_r(x-y),
\quad r=1,2
,3,4.
$$
Changing
here $x\leftrightarrow y$, one gets
$[rr11]'-[rr11]= \theta_1(u+v)\theta_1(u-v)\theta_r(x+y)\theta_r(x-y)
=[11rr]-[11rr]'$.
Similarly, (\ref{si1d})--(\ref{si1f}) follow from (\ref{tt2c})--(\ref{tt2a}), respectively.
To prove (\ref{si2a}), we write (\ref{si1a}) in terms of
the variables $u,v,x,y$ and
exchange $u\leftrightarrow x$. Then
$$[2211]-[1122]=
-\theta_1(u+v)\theta_1(u-v)\theta_2(x+y)\theta_2(x-y)+
\theta_2(u+v)\theta_2(u-v)\theta_1(x+y)\theta_1(x-y).$$
Now (\ref{si2a}) holds by virtue of
(\ref{tt2a}). Identities (\ref{si2b})--(\ref{si2f}) can be proved
in a similar way.
Finally, it is easy to see that
(\ref{si4}) follows from (\ref{tt5}).
Indeed, subtracting (\ref{tt5}) from the same identity with
the exchange ${x\leftrightarrow y}$ yields (\ref{ft11a}).
All other identities in (\ref{si4})
are proved in a similar way.
}
\begin{rem}
Identities (\ref{si1}), (\ref{si2})
differ slightly from those written by Jacobi.
For example, in {\rm\cite[p. 507]{J}} one can find the relations
$[1122]+[4433]=[2211]'+[3344]'$,
$[1122]-[4433]=[1122]'-[4433]'$
which are appropriate linear combinations of (\ref{si1a}), (\ref{si1f}), (\ref{si2a}), (\ref{si2f}).
We should also stress that these Jacobi identities are direct corollaries of
(\ref{tt1}) at $r=2$ and (\ref{tt2a}).
This is in complete agreement with the derivation
of (\ref{si1}), (\ref{si2}) from the Weierstrass addition formulae.
\end{rem}

\subsection{Five term Riemann identities}
The Riemann identities (the term is due to Mumford
\cite[page 20]{M1}) are simple corollaries of the Jacobi relations
(\ref{ft2}), (\ref{si1})--(\ref{si4}).
They each express a ``primed'' quantity as a linear combination of
some appropriate four
``unprimed'' ones.

\smallskip\noindent
Hence from (\ref{ft2}) we have the four simplest Riemann identities:
\bs{j1}
\ba
2[1]'&=[1]+[2]-[3]+[4],\label{j1a}\\
2[2]'&=[1]+[2]+[3]-[4],\label{j1b}\\
\hspace{-4cm} {\bf R.I} \hspace{2cm}
2[3]'&=-[1]+[2]+[3]+[4],\label{j1c}\\
2[4]'&=[1]-[2]+[3]+[4].\label{j1d}
\end{align}
\es
Let us emphasize that (\ref{j1}) is equivalent to (\ref{ft2}).
\smallskip

\noindent
Next we list all possible (twelve) identities that are
obtained from
(\ref{si1}), (\ref{si2}):
\bs{j2}
\ba
2[1122]'=[1122]+[2211]+[3344]-[4433],\label{j2a}\\
2[1133]'=[1133]+[3311]+[2244]-[4422],\label{j2b}\\
2[1144]'=[1144]+[4411]+[2233]-[3322].\label{j2c}
\end{align}
\es
\vspace{-0.5cm}
\bs{j3}
\ba
2[2211]'=[2211]+[1122]+[4433]-[3344],\label{j3a}\\
2[2233]'=[2233]+[3322]+[1144]-[4411],\label{j3b}\\
\hspace{-4cm} {\bf R.II} \hspace{2cm}
2[2244]'=[2244]+[4422]+[1133]-[3311].\label{j3c}
\end{align}
\es
\vspace{-0.5cm}
\bs{j4}
\ba
2[3311]'=[3311]+[1133]+[4422]-[2244],\label{j4a}\\
2[3322]'=[3322]+[2233]+[4411]-[1144],\label{j4b}\\
2[3344]'=[3344]+[4433]+[1122]-[2211].\label{j4c}
\end{align}
\es
\bs{j5}
\vspace{-0.5cm}
\ba
2[4411]'=[4411]+[1144]+[3322]-[2233],\label{j5a}\\
2[4422]'=[4422]+[2244]+[3311]-[1133],\label{j5b}\\
2[4433]'=[4433]+[3344]+[2211]-[1122].\label{j5c}
\end{align}
\es
\smallskip

\noindent
Clearly, the system (\ref{j2})--(\ref{j5}) is equivalent to (\ref{si1}),
(\ref{si2}). Finally, identities (\ref{si4})
are equivalent to
\bs{j6}
\ba
2[1234]'=-[1234]+[2143]+[3421]+[4312],\label{j6a}\\
2[2143]'=-[2143]+[1234]+[3412]+[4321],\label{j6b}\\
\hspace{-4cm} {\bf R.III} \hspace{2cm}
2[3412]'=-[3412]+[4321]+[1234]+[2143],\label{j6c}\\
2[4321]'=-[4321]+[3412]+[1234]+[2143].\label{j6d}
\end{align}
\es
\begin{rem}
The identities presented here essentially coincide with the ones
given by Mumford {\rm\cite[p. 20]{M1}}. See also {\rm \cite{WW}}.
\end{rem}

\subsection{Equivalence of addition formulae and Jacobi identities}
In section \ref{qr4}, we have obtained the Jacobi identities from
the addition formulae.
In its turn, one can show that the system (\ref{ft2}),
(\ref{si1}), (\ref{si2}) implies
the addition formulae (\ref{tt2})-(\ref{tt5}). The proof is similar to the
one given by Koornwinder for the Riemann identities \cite{K}.

\smallskip

{\small
\noindent
In accordance with (\ref{ft'}), the
relation (\ref{ft2a}) acquires the form
\bel{z1}
\hspace{-1.2cm}
\theta_1(u+x)\theta_1(u-x)\theta_1(v+y)\theta_1(v-y)+
\theta_2(u+x)\theta_2(u-x)\theta_2(v+y)\theta_2(v-y)\\ \hspace{-1.2cm}
=\, \theta_1(v+x)\theta_1(v-x)\theta_1(u+y)\theta_1(u-y)+
\theta_2(v+x)\theta_2(v-x)\theta_2(u+y)\theta_2(u-y).
\ee
Changing here $u\leftrightarrow x$ and
$v\leftrightarrow x$, one obtains two
additional relations:
\bel{z2}
\hspace{-1.3cm}
-\theta_1(u+x)\theta_1(u-x)\theta_1(v+y)\theta_1(v-y)+
\theta_2(u+x)\theta_2(u-x)\theta_2(v+y)\theta_2(v-y)\\ \hspace{-1.4cm}
=-\theta_1(u+v)\theta_1(u-v)\theta_1(x+y)\theta_1(x-y)+
\theta_2(u+v)\theta_2(u-v)\theta_2(x+y)\theta_2(x-y),
\ee
\bel{z3}
\hspace{-1cm}
\theta_1(u+v)\theta_1(u-v)\theta_1(x+y)\theta_1(x-y)+
\theta_2(u+v)\theta_2(u-v)\theta_2(x+y)\theta_2(x-y)\\ \hspace{-1.4cm}
=-\theta_1(v+x)\theta_1(v-x)\theta_1(u+y)\theta_1(u-y)+
\theta_2(v+x)\theta_2(v-x)\theta_2(u+y)\theta_2(u-y).
\ee
Introduce the notation:
\bel{z4}
A_j:=\theta_j(u+x)\theta_j(u-x)\theta_j(v+y)\theta_j(v-y),\\
B_j:=\theta_j(u+y)\theta_j(u-y)\theta_j(v+x)\theta_j(v-x),\\
C_j:=\theta_j(u+v)\theta_j(u-v)\theta_j(x+y)\theta_j(x-y).
\ee
Then relations (\ref{z1})--(\ref{z3})
acquire the form $A_1-B_1=B_2-A_2$,
$A_1-C_1=A_2-C_2$, $B_1+C_1=B_2-C_2$ which is
a system
of linear equations for the unknowns $A_2, B_2, C_2$.
The system is degenerate with
compatibility condition $A_1-B_1=C_1$.
In terms of the theta functions, this condition is nothing but
equation (\ref{tt3})
with $r=1$. Equivalently,
one can treat the above equations for $A_j,B_j,C_j$
as a linear system for the unknowns
$A_1, B_1, C_1$. Then, for example, $C_1=B_2-A_2$ which is (\ref{tt3})
with $r=2$. The
other addition formulae can be obtained in a similar way.
}

\section{Particular identities}
\subsection{Consequences of the bilinear identities}
One can obtain twelve
particular identities from the general system (\ref{es1}) putting $v=0$ or
$v=\pm u$ (actually, the restriction $v=-u$
can be applied only for identity
(\ref{es1d}) which leads to (\ref{bc3c}) below):
\bs{bc1}
\ba
&2\theta_1^2(u|2\tau)=\theta_4(u|\tau)\theta_3(0|\tau)
-\theta_3(u|\tau)\theta_4(0|\tau),\label{bc1a}\\
&2\theta_2^2(u|2\tau)=\theta_3(u|\tau)\theta_3(0|\tau)
-\theta_4(u|\tau)\theta_4(0|\tau),\label{bc1b}\\
&2\theta_3^2(u|2\tau)=\theta_3(u|\tau)\theta_3(0|\tau)
+\theta_4(u|\tau)\theta_4(0|\tau),\label{bc1c}\\
&2\theta_4^2(u|2\tau)=\theta_3(u|\tau)\theta_4(0|\tau)
+\theta_4(u|\tau)\theta_3(0|\tau),\label{bc1d}
\end{align}
\es
\vspace{-0.6cm}
\bs{bc2}
\ba
&2\theta_1(u|2\tau)\theta_4(u|2\tau)=
\theta_1(u|\tau)\theta_2(0|\tau),\label{bc2a}\\
&2\theta_2(u|2\tau)\theta_3(u|2\tau)=
\theta_2(u|\tau)\theta_2(0|\tau),\label{bc2b}
\end{align}
\es
\vspace{-0.6cm}
\bs{bc3}
\ba
&2\theta_2(2u|2\tau)\theta_2(0|2\tau)=\theta_3^2(u|\tau)
-\theta_4^2(u|\tau),\label{bc3a}\\
&2\theta_2(2u|2\tau)\theta_3(0|2\tau)=
\theta_2^2(u|\tau)-\theta_1^2(u|\tau),\label{bc3b}\\
&2\theta_3(2u|2\tau)\theta_2(0|2\tau)=
\theta_2^2(u|\tau)+\theta_1^2(u|\tau),\label{bc3c}\\
&2\theta_3(2u|2\tau)\theta_3(0|2\tau)=\theta_3^2(u|\tau)+
\theta_4^2(u|\tau),\label{bc3d}
\end{align}
\es
\vspace{-0.6cm}
\bs{bc4}
\ba
&\,\,\, \theta_1(2u|2\tau)\theta_4(0|2\tau)=
\theta_1(u|\tau)\theta_2(u|\tau),\label{bc4a}\\
&\,\,\, \theta_4(2u|2\tau)\theta_4(0|2\tau)=
\theta_3(u|\tau)\theta_4(u|\tau).\label{bc4b}
\end{align}
\es

\noindent
In \cite[Section 21.52]{WW}
there are two particular equations relating
theta functions
with modular parameters $\tau$ and $2\tau$ which are
called {\it the transformations of
Landen's type}:
\bse\label{lt1}
\bel{lt1a}
\frac{\theta_4(2u|2\tau)}{\theta_4(0|2\tau)}=\frac{\theta_3(u|\tau)
\theta_4(u|\tau)}{\theta_3(0|\tau)\theta_4(0|\tau)},
\ee
\vspace{-0.3cm}
\bel{lt1b}
\frac{\theta_1(2u|2\tau)}{\theta_4(0|2\tau)}=\frac{\theta_1(u|\tau)
\theta_2(u|\tau)}{\theta_3(0|\tau)\theta_4(0|\tau)}.
\ee
\ese
The first identity is derived from (\ref{bc4b}) and from the relation
\bel{lt2}
\theta_4^2(0|2\tau)=\theta_3(0|\tau)\theta_4(0|\tau)
\ee
which is also a corollary of (\ref{bc4b}).
The identity (\ref{lt1b}) is
a ratio of (\ref{bc4a}) and (\ref{lt2}).\\
Note also that (\ref{lt1b}) can be obtained from (\ref{lt1a}) by the shift
$u\to u+\frac{\tau}{2}$.

\subsection{Particular addition formulae}\label{paf}

One can obtain important particular cases of
the Weierstrass addition formulae
which include two variables.
Here we present the complete list of
eighteen identities which easily follow from (\ref{tt1})--(\ref{tt5}):
\bse\label{ad1}
\bel{ad1a}
\theta_1(u+v)\theta_1(u-v)\theta_2^2(0) & =\,\,
\theta_1^2(u)\theta_2^2(v)-\theta_2^2(u)\theta_1^2(v)\\
&=\,\, \theta_4^2(u)\theta_3^2(v)-\theta_3^2(u)\theta_4^2(v),
\ee
\bel{ad1b}
\theta_1(u+v)\theta_1(u-v)\theta_3^2(0) & =\,\,
\theta_1^2(u)\theta_3^2(v)-\theta_3^2(u)\theta_1^2(v)\\
& =\,\, \theta_4^2(u)\theta_2^2(v)-\theta_2^2(u)\theta_4^2(v),
\ee
\bel{ad1c}
\theta_1(u+v)\theta_1(u-v)\theta_4^2(0) & =\,\,
\theta_1^2(u)\theta_4^2(v)-\theta_4^2(u)\theta_1^2(v)\\
& =\,\, \theta_3^2(u)\theta_2^2(v)-\theta_2^2(u)\theta_3^2(v),
\ee
\ese
\bse\label{ad2}
\bel{ad2a}
\theta_2(u+v)\theta_2(u-v)\theta_2^2(0)& = \,\,
\theta_2^2(u)\theta_2^2(v)-\theta_1^2(u)\theta_1^2(v)\\
& =\,\, \theta_3^2(u)\theta_3^2(v)-\theta_4^2(u)\theta_4^2(v),
\ee
\be
\theta_2(u+v)\theta_2(u-v)\theta_3^2(0) & =\,\,
\theta_3^2(u)\theta_2^2(v)-\theta_1^2(u)\theta_4^2(v)\\
&=\,\, \theta_2^2(u)\theta_3^2(v)-\theta_4^2(u)\theta_1^2(v),
\ee
\be
\theta_2(u+v)\theta_2(u-v)\theta_4^2(0) & =\,\,
\theta_4^2(u)\theta_2^2(v)-\theta_1^2(u)\theta_3^2(v)\\
&=\,\, \theta_2^2(u)\theta_4^2(v)-\theta_3^2(u)\theta_1^2(v),
\ee
\ese
\bse\label{ad3}
\bel{ad3a}
\theta_3(u+v)\theta_3(u-v)\theta_2^2(0)& =\,\,
\theta_2^2(u)\theta_3^2(v)+\theta_1^2(u)\theta_4^2(v)\\
& = \,\,\theta_3^2(u)\theta_2^2(v)+\theta_4^2(u)\theta_1^2(v),
\ee
\bel{ad3b}
\theta_3(u+v)\theta_3(u-v)\theta_3^2(0) & =\,\,
\theta_1^2(u)\theta_1^2(v)+\theta_3^2(u)\theta_3^2(v)\\
& =\,\, \theta_2^2(u)\theta_2^2(v)+\theta_4^2(u)\theta_4^2(v),
\ee
\bel{ad3c}
\theta_3(u+v)\theta_3(u-v)\theta_4^2(0) & = \,\,
\theta_4^2(u)\theta_3^2(v)-\theta_1^2(u)\theta_2^2(v)\\
&=\,\, \theta_3^2(u)\theta_4^2(v)-\theta_2^2(u)\theta_1^2(v),
\ee
\ese
\bse\label{ad4}
\bel{ad4a}
\theta_4(u+v)\theta_4(u-v)\theta_2^2(0)& =\,\,
\theta_1^2(u)\theta_3^2(v)+\theta_2^2(u)\theta_4^2(v)\\
& = \,\, \theta_3^2(u)\theta_1^2(v)+\theta_4^2(u)\theta_2^2(v),
\ee
\bel{ad4b}
\theta_4(u+v)\theta_4(u-v)\theta_3^2(0)& =
\theta_1^2(u)\theta_2^2(v)+\theta_3^2(u)\theta_4^2(v)\\
& = \,\, \theta_2^2(u)\theta_1^2(v)+\theta_4^2(u)\theta_3^2(v),
\ee
\bel{ad4c}
\theta_4(u+v)\theta_4(u-v)\theta_4^2(0)& =\,\,
\theta_4^2(u)\theta_4^2(v)-\theta_1^2(u)\theta_1^2(v)\\
& = \,\, \theta_3^2(u)\theta_3^2(v)-\theta_2^2(u)\theta_2^2(v),
\ee
\ese
\bs{ad5}
\ba
\theta_1(u+v)\theta_2(u-v)\theta_3(0)\theta_4(0)=
\theta_1(u)\theta_2(u)\theta_3(v)\theta_4(v)+
\theta_3(u)\theta_4(u)\theta_1(v)\theta_2(v),\\
\theta_1(u+v)\theta_3(u-v)\theta_2(0)\theta_4(0)=
\theta_1(u)\theta_3(u)\theta_2(v)\theta_4(v)+
\theta_2(u)\theta_4(u)\theta_1(v)\theta_3(v),\\
\theta_1(u+v)\theta_4(u-v)\theta_2(0)\theta_3(0)=
\theta_1(u)\theta_4(u)\theta_2(v)\theta_3(v)+
\theta_2(u)\theta_3(u)\theta_1(v)\theta_4(v),\\
\theta_2(u+v)\theta_3(u-v)\theta_2(0)\theta_3(0)=
\theta_2(u)\theta_3(u)\theta_2(v)\theta_3(v)-
\theta_1(u)\theta_4(u)\theta_1(v)\theta_4(v),\\
\theta_2(u+v)\theta_4(u-v)\theta_2(0)\theta_4(0)=
\theta_2(u)\theta_4(u)\theta_2(v)\theta_4(v)-
\theta_1(u)\theta_3(u)\theta_1(v)\theta_3(v),\\
\theta_3(u+v)\theta_4(u-v)\theta_3(0)\theta_4(0)=
\theta_3(u)\theta_4(u)\theta_3(v)\theta_4(v)-
\theta_1(u)\theta_2(u)\theta_1(v)\theta_2(v).
\end{align}
\es
\begin{rem}
The complete list of identities (\ref{ad1})--(\ref{ad5})
was originally obtained by Jacobi
{\rm \cite[p. 510] {J}} as a particular
specification  of identities (\ref{ft2}),
(\ref{si1})--(\ref{si4}),
see also {\rm\cite[pp. 76-78] {W}, \cite[487-488]{WW}}.
Mumford {\rm\cite[p. 22]{M1}}
has obtained a part of relations (\ref{ad1})-(\ref{ad5}) as
specific cases of
the Riemann identities
(\ref{j1})--(\ref{j6}).
\end{rem}

\noindent
As a byproduct of (\ref{ad1})--(\ref{ad4}), one gets some extra identities:
\bs{ad6}
\ba
\theta_1^2(u)\theta_1^2(v)-\theta_2^2(u)\theta_2^2(v)=
\theta_4^2(u)\theta_4^2(v)-\theta_3^2(u)\theta_3^2(v),\\
\theta_1^2(u)\theta_2^2(v)-\theta_2^2(u)\theta_1^2(v)=
\theta_4^2(u)\theta_3^2(v)-\theta_3^2(u)\theta_4^2(v),\\
\theta_1^2(u)\theta_3^2(v)-\theta_3^2(u)\theta_1^2(v)=
\theta_4^2(u)\theta_2^2(v)-\theta_2^2(u)\theta_4^2(v),\\
\theta_1^2(u)\theta_4^2(v)-\theta_4^2(u)\theta_1^2(v) =
\theta_3^2(u)\theta_2^2(v)-\theta_2^2(u)\theta_3^2(v).
\end{align}
\es
In particular, the following identity holds:
\bel{ad7}
\theta_1^4(u)+\theta_3^4(u)=\theta_2^4(u)+\theta_4^4(u).
\ee
Certainly, the last relation is a corollary of (\ref{bc3}).
Thus, in addition to (\ref{tc1}), one gets
another famous identity for theta constants:
\bel{tc2}
\theta_3^4(0)=\theta_2^4(0)+\theta_4^4(0).
\ee

\subsection{Duplication formulae}
\label{sdf}
The duplication formulae relate the
functions $\theta_a(2u|\tau),\,a=1,2,3,4$ with
appropriate combinations of the functions $\theta_b(u|\tau)$. All these identities emerge as
further degenerations of addition formulae (\ref{ad1})--(\ref{ad5}).
Here is the complete list:
\bel{df1}
\theta_1(2u)\theta_2(0)\theta_3(0)\theta_4(0)=2\,
\theta_1(u)\theta_2(u)\theta_3(u)\theta_4(u),
\ee
\vspace{-0.4cm}
\bs{df2}
\begin{align}
\theta_2(2u)\theta_2(0)\theta_3^2(0)&=\,
\theta_2^2(u)\theta_3^2(u)-\theta_1^2(u)\theta_4^2(u),\label{df2b}\\
\theta_2(2u)\theta_2(0)\theta_4^2(0)&=\,
\theta_2^2(u)\theta_4^2(u)-\theta_1^2(u)\theta_3^2(u),\label{df2a}\\
\theta_2(2u)\theta_2^3(0)&= \,\theta_2^4(u)-\theta_1^4(u),\label{df2c}\\
\theta_2(2u)\theta_2^3(0)&= \, \theta_3^4(u)-\theta_4^4(u)\label{df2d}.
\end{align}
\es
\vspace{-0.4cm}
\bs{df3}
\begin{align}
\theta_3(2u)\theta_3(0)\theta_2^2(0)&=\,
\theta_2^2(u)\theta_3^2(u)+\theta_1^2(u)\theta_4^2(u),\label{df3b}\\
\theta_3(2u)\theta_3(0)\theta_4^2(0)&=\,
\theta_3^2(u)\theta_4^2(u)-\theta_1^2(u)\theta_2^2(u),\label{df3a}\\
\theta_3(2u)\theta_3^3(0)&=\, \theta_1^4(u)+\theta_3^4(u),\label{df3c}\\
\theta_3(2u)\theta_3^3(0)&=\, \theta_2^4(u)+\theta_4^4(u)\label{df3d}.
\end{align}
\es
\vspace{-0.4cm}
\bs{df4}
\begin{align}
\theta_4(2u)\theta_4(0)\theta_2^2(0)&=\,
\theta_2^2(u)\theta_4^2(u)+\theta_1^2(u)\theta_3^2(u),\label{df4b}\\
\theta_4(2u)\theta_4(0)\theta_3^2(0)&=\,
\theta_1^2(u)\theta_2^2(u)+\theta_3^2(u)\theta_4^2(u),\label{df4a}\\
\theta_4(2u)\theta_4^3(0)&=\, \theta_4^4(u)-\theta_1^4(u),\label{df4c}\\
\theta_4(2u)\theta_4^3(0)&=\, \theta_3^4(u)-\theta_2^4(u)\label{df4d}.
\end{align}
\es
One can unify
the sets of equations \{(\ref{df2b}), (\ref{df3a}), (\ref{df4b})\};
\{(\ref{df2a}), (\ref{df3b}),(\ref{df4a})\};
\{(\ref{df2c}), (\ref{df3c}),(\ref{df4c})\} and
\{(\ref{df2d}), (\ref{df3d}),(\ref{df4d})\} by writing down
all the 12 identities
(\ref{df2})--(\ref{df4}) in the compressed form
\bs{df5}
\begin{align}
(-1)^{\beta+\gamma}\theta_1^2(u)\theta_{\al+1}^2(u)+
\theta_{\beta+1}^2(u)\theta_{\gamma+1}^2(u) &=\,
\theta_{\beta+1}(2u)\theta_{\beta+1}(0)\theta_{\gamma+1}^2(0),
\label{df5a}\\
(-1)^{\beta+\gamma}\theta_1^2(u)\theta_{\al+1}^2(u)-
\theta_{\beta+1}^2(u)\theta_{\gamma+1}^2(u) &=-\,
\theta_{\gamma+1}(2u)\theta_{\gamma+1}(0)\theta_{\beta+1}^2(0),
\label{df5b}\\
\theta_{\al+1}(2u)\theta^3_{\al+1}(0) &=\, \theta^4_{\al+1}(u)+(-1)^\al
\theta_1^4(u)\,,\label{df5c}\\
\theta_{\al+1}(2u)\theta^3_{\al+1}(0) &=\,
(-1)^{\gamma+1}\theta^4_{\beta+1}(u)+
(-1)^{\beta+1}\theta_{\gamma+1}^4(u).\label{df5d}
\end{align}
\es
where in (\ref{df5a}), (\ref{df5b}), and (\ref{df5d})
the indices $\al,\beta,\gamma$ are assumed to be any cyclic permutation
of $\{1,2,3\}$ and in (\ref{df5c}) $\al=1,2,3$. Together
with  (\ref{df1}), the
system (\ref{df5}) yields the complete set of duplication formulae.

\bibliographystyle{12}
\bibliographystyle{amsalpha}

\end{document}